\documentclass[12pt,reqno]{amsart}
\usepackage{amsfonts,amsmath,amssymb}
\usepackage{graphicx}
\usepackage{psfrag}
\usepackage{pstricks}
\usepackage{pst-tree}

\allowdisplaybreaks
\advance\textheight by \topskip

\newtheorem{theorem}{Theorem}

\numberwithin{equation}{section}

\title[A new (?) continued fraction expansion for the reciprocal of a 
	$q$-series]{A new (?) continued fraction expansion for the reciprocal of a 
	$\boldsymbol q$-series}

\author[H.~Prodinger]{Helmut Prodinger}
\address{Helmut Prodinger\\
Department of Mathematics\\
University of Stellenbosch\\
7602 Stellenbosch\\
South Africa}
\email{hproding@sun.ac.za}

\date{June 4, 2008}

\begin{document}
\begin{abstract}
We prove a continued fraction expansion for the reciprocal of a certain $q$-series.
All the specialists in the world are asked whether it is new or not.
\end{abstract}

\maketitle

\section{Introduction}

Let
\begin{equation*}
(x;q)_n=(1-x)(1-xq)\dots(1-xq^{n-1}).
\end{equation*}
Series of the type
\begin{equation*}
\sum \frac{(a_1;q)_n\dots(a_s;q)_n}{(b_1;q)_n\dots(b_t;q)_n}z^n
\end{equation*}
are called $q$-series (basic hypergeometric series), and the summation might be either over the nonnegative integers or over
the whole set of integers. 

There are several results about such series and \emph{continued fraction expansions};
to give one example, set
\begin{equation*}
_2\psi_2(\alpha, \beta;\gamma,\delta\mid z)=\sum_{n\in\mathbb{Z}}
\frac{(\alpha;q)_n(\beta;q)_n}{(\gamma;q)_n(\delta;a)_n}z^n,
\end{equation*}
then 
\begin{equation*}
\frac{_2\psi_2(\alpha, \beta;\gamma,\delta\mid z)}{_2\psi_2(\alpha, \beta q;\gamma,\delta q\mid z)}
\end{equation*}
is expanded into a continued fraction in \cite{Denis90}.

The examples that I could trace in the literature usually deal with \emph{quotients} of \emph{similar} functions.  The is a school of indian mathematicians who deal with such questions; \textsf{MathSciNet} will easily produce a list of such papers (which I found hard to get!).

The continued fraction that I will present here only deals with \emph{one} $q$-series, not a quotient. I have not been able to find it in the literature, although it might exist.

Thus, I post this draft now, and ask all the specialists in the world for feedback. If nobody will
identify it within a month, say, then I will dare to submit it to a journal. Otherwise, I will label it
\emph{written for my own benefit}.

\section{Main result}

\begin{theorem} Let 
\begin{equation*}
G(z)=\sum_{n\ge0}\frac{(y;q)_nz^{n}}{(x;q)_n}.
\end{equation*}
Then we have the continued fraction expansion
\begin{equation*}
\frac{z}{G(z)}
=
\cfrac{z}{a_1+
	\cfrac{z}{a_2+
		\cfrac{z}{a_3+
			\cfrac{z}{a_4+\cfrac z{\ddots}
			}
		}
	}
}
\end{equation*}
with $a_1=1$, and for $k\ge1$
\begin{align*}
a_{2k}&=\frac{(x;q)_{k-1}(1-xq^{2k-2})(yq;q)_{k-1}}{(1-y)(yq)^{k-1}(\frac xy;q)_{k-1}(q;q)_{k-1}},\\
a_{2k+1}&=-\frac{(1-y)y^{k-1}(\frac xy;q)_{k-1}(1-xq^{2k-1})(q;q)_{k-1}}{(x;q)_{k}(yq;q)_{k}}.
\end{align*}
\end{theorem}
\emph{Proof.}
We write
\begin{align*}
\frac {z}{G(z)}&=\cfrac{z}{N_0}=\cfrac{z}{a_1+\cfrac{z}{N_1}}=\cfrac{z}{a_1+\cfrac{z}{a_2+\cfrac{z}{N_2}}}=\dots,
\end{align*}
and set
\begin{equation*}
N_i=\frac{r_i}{s_i}.
\end{equation*}
This means
\begin{equation*}
N_i=a_{i+1}+\frac{z}{N_{i+1}}
\end{equation*}
or
\begin{equation*}
\frac{z}{N_{i+1}}=\frac{zs_{i+1}}{r_{i+1}}=N_i-a_{i+1}=\frac{r_i}{s_i}-a_{i+1}
=\frac{r_i-a_{i+1}s_i}{s_i}.
\end{equation*}
We may identify numerators and denominators and deduce the recursion
\begin{equation*}
zs_{i+1}=s_{i-1}-a_{i+1}s_i.
\end{equation*}
Since $a_1$ is an exceptional value, we only start the recursion with
\begin{equation*}
s_0=1, \qquad s_1=\sum_{n\ge0}\frac{(y;q)_{n+1}z^{n}}{(x;q)_{n+1}}.
\end{equation*}
Note that the numbers $a_i$ are uniquely determined by annihilating the constant term in
$s_{i-1}-a_{i+1}s_i$, making $s_{i+1}$ a power series expansion. Our claim follows now
by the following \emph{explicit} formul\ae{} (for $k\ge0$)
\begin{align*}
s_{2k}&=\frac{(-1)^kq^{ \binom k2}}{(q;q)_{k-1}}\sum_{n\ge0}\frac{z^n(q^{n+1};q)_{k-1}(yq;q)_{n+k}}{(xq^k;q)_{n+k}},\\*
s_{2k+1}&=(1-y)y^k(\tfrac xy;q)_k(-1)^kq^{\binom{k+1}2}\sum_{n\ge0}\frac{z^n(q^{n+1};q)_{k}(yq^{k+1};q)_n}{(x;q)_{n+2k+1}},
\end{align*}
provided we are able to establish these formul\ae{} by induction via the recursion. 
The initial values follow by inspection, and the induction step must be split into two
computations, according to the parity of the indices.
\small
\begin{align*}
s_{2k}&-a_{2k+2}s_{2k+1}=\frac{(-1)^kq^{ \binom k2}}{(q;q)_{k-1}}\sum_{n\ge0}\frac{z^n(q^{n+1};q)_{k-1}(yq;q)_{n+k}}{(xq^k;q)_{n+k}}\\
&-\frac{(x;q)_{k}(qy;q)_{k}(1-xq^{2k})}{(qy)^{k}(\frac xy;q)_{k}(1-y)(q;q)_{k}}(1-y)y^k(\tfrac xy;q)_k(-1)^kq^{\binom{k+1}2}\sum_{n\ge0}\frac{z^n(q^{n+1};q)_{k}(yq^{k+1};q)_n}{(x;q)_{n+2k+1}}\\
&=\frac{(-1)^kq^{ \binom k2}}{(q;q)_{k-1}}\sum_{n\ge0}\frac{z^n(q^{n+1};q)_{k-1}(yq;q)_{n+k}}{(xq^k;q)_{n+k}}-\frac{(1-xq^{2k})}{(q;q)_{k}}(-1)^kq^{\binom{k}2}\sum_{n\ge0}\frac{z^n(q^{n+1};q)_{k}(yq;q)_{n+k}}{(xq^k;q)_{n+k+1}}\\
&=\frac{(-1)^kq^{ \binom k2}}{(q;q)_{k}}\bigg[
(1-q^k)\sum_{n\ge0}\frac{z^n(q^{n+1};q)_{k-1}(yq;q)_{n+k}}{(xq^k;q)_{n+k}}-(1-xq^{2k})\sum_{n\ge0}\frac{z^n(q^{n+1};q)_{k}(yq;q)_{n+k}}{(xq^k;q)_{n+k+1}}\bigg]\\
&=\frac{(-1)^kq^{ \binom k2}}{(q;q)_{k}}
\sum_{n\ge0}\frac{z^n(q^{n+1};q)_{k-1}(yq;q)_{n+k}}{(xq^k;q)_{n+k+1}}\Big[(1-q^k)(1-xq^{n+2k})-(1-xq^{2k})(1-q^{n+k})\Big]\\
&=\frac{(-1)^kq^{ \binom k2}}{(q;q)_{k}}
\sum_{n\ge0}\frac{z^n(q^{n+1};q)_{k-1}(yq;q)_{n+k}}{(xq^k;q)_{n+k+1}}(-1)q^k(1-xq^k)(1-q^n)\\
&=\frac{(-1)^{k+1}q^{ \binom {k+1}2}}{(q;q)_{k}}
\sum_{n\ge0}\frac{z^n(q^{n};q)_{k}(yq;q)_{n+k}}{(xq^{k+1};q)_{n+k}}\\
&=zs_{2k+2};
\end{align*}
\begin{align*}
s_{2k-1}&-a_{2k+1}s_{2k}=
(1-y)y^{k-1}(\tfrac xy;q)_{k-1}(-1)^{k-1}q^{\binom{k}2}\sum_{n\ge0}\frac{z^n(q^{n+1};q)_{k-1}(yq^{k};q)_n}{(x;q)_{n+2k-1}}\\
&\quad+\frac{(1-y)y^{k-1}(\frac xy;q)_{k-1}(q;q)_{k-1}(1-xq^{2k-1})}{(qy;q)_{k}(x;q)_{k}}\frac{(-1)^kq^{ \binom k2}}{(q;q)_{k-1}}\sum_{n\ge0}\frac{z^n(q^{n+1};q)_{k-1}(yq;q)_{n+k}}{(xq^k;q)_{n+k}}\\
&=
(1-y)y^{k-1}(\tfrac xy;q)_{k-1}(-1)^{k-1}q^{\binom{k}2}\sum_{n\ge0}\frac{z^n(q^{n+1};q)_{k-1}(yq^{k};q)_n}{(x;q)_{n+2k-1}}\\
&\quad+(1-y)y^{k-1}(\tfrac xy;q)_{k-1}(1-xq^{2k-1})(-1)^kq^{ \binom k2}\sum_{n\ge0}\frac{z^n(q^{n+1};q)_{k-1}(yq^{k+1};q)_{n}}{(x;q)_{n+2k}}\\
&=(1-y)y^{k-1}(\tfrac xy;q)_{k-1}(-1)^{k-1}q^{\binom{k}2}\\&\quad\times\bigg[\sum_{n\ge0}\frac{z^n(q^{n+1};q)_{k-1}(yq^{k};q)_n}{(x;q)_{n+2k-1}}-(1-xq^{2k-1})\sum_{n\ge0}\frac{z^n(q^{n+1};q)_{k-1}(yq^{k+1};q)_{n}}{(x;q)_{n+2k}}\bigg]\\
&=(1-y)y^{k-1}(\tfrac xy;q)_{k-1}(-1)^{k-1}q^{\binom{k}2}\sum_{n\ge0}\frac{z^n(q^{n+1};q)_{k-1}(yq^{k+1};q)_{n-1}}{(x;q)_{n+2k}}\\
&\quad\times\Big[(1-xq^{n+2k-1})(1-yq^k)-(1-xq^{2k-1})(1-yq^{k+n})\Big]\\
&=(1-y)y^{k-1}(\tfrac xy;q)_{k-1}(-1)^{k-1}q^{\binom{k}2}\sum_{n\ge0}\frac{z^n(q^{n+1};q)_{k-1}(yq^{k+1};q)_{n-1}}{(x;q)_{n+2k}}\\&\quad\times(-1)q^{k-1}(1-q^n)yq(1-\tfrac xyq^{k-1})\\
&=(1-y)y^{k}(\tfrac xy;q)_{k}(-1)^{k}q^{\binom{k+1}2}\sum_{n\ge0}\frac{z^n(q^{n};q)_{k}(yq^{k+1};q)_{n-1}}{(x;q)_{n+2k}}\\
&=zs_{2k+1}.
\end{align*}
\normalsize

\emph{Remark.} Computer experiments indicate that we cannot add additional factors $(u;q)_n$
etc. in either numerator or denominator, as then the expressions for $a_i$ become very messy and
don't factor nicely.

\medskip

We note two special cases explicitly. Set $y=0$, then
\begin{align*}
a_{2k}&=\frac{(x;q)_{k-1}(1-xq^{2k-2})}{x^{k-1}q^{\binom k2}(q;q)_{k-1}},\\
a_{2k+1}&=-\frac{x^{k-1}q^{\binom {k-1}2}(1-xq^{2k-1})(q;q)_{k-1}}{(x;q)_{k}}.
\end{align*}
Set $x=0$, then
\begin{align*}
a_{2k}&=\frac{(yq;q)_{k-1}}{(1-y)(yq)^{k-1}(q;q)_{k-1}},\\
a_{2k+1}&=-\frac{(1-y)y^{k-1}(q;q)_{k-1}}{(yq;q)_{k}}.
\end{align*}

\bibliographystyle{amsplain} 



\providecommand{\bysame}{\leavevmode\hbox to3em{\hrulefill}\thinspace}
\providecommand{\MR}{\relax\ifhmode\unskip\space\fi MR }
\providecommand{\MRhref}[2]{%
  \href{http://www.ams.org/mathscinet-getitem?mr=#1}{#2}
}
\providecommand{\href}[2]{#2}

\end{document}